\documentclass[12pt]{amsart}
\usepackage{amsthm, amsmath}
\usepackage{amssymb} 
\usepackage{graphicx}
\usepackage{epsf}
\usepackage[all]{xy}

\theoremstyle{plain}
\newtheorem{theorem}{Theorem}[section]
\newtheorem{lemma}[theorem]{Lemma}
\newtheorem{cor}[theorem]{Corollary}
\newtheorem{prop}[theorem]{Proposition}

\theoremstyle{remark}

\theoremstyle{definition}

\newcommand{\CC}{\mathbb C}
\newcommand{\ZZ}{\mathbb Z}
\newcommand{\QQ}{\mathbb Q}

\author{Thomas J. Haines}
\title{Structure constants for Hecke\\
and representation rings}

\date{}

\begin{document}

\maketitle

\begin{abstract}
In \cite{KLM} the authors study certain structure constants for two related rings: the spherical Hecke algebra of a split connected reductive group over a local non-Archimedean field, and the representation ring of the Langlands dual group.  The former are defined relative to characteristic functions of double cosets, and the latter relative to highest weight representations.  They prove that the nonvanishing of one of the latter structure constants always implies the nonvanishing of the corresponding former one.  For ${\rm GL}_n$, the reverse implication also holds, and is due to P. Hall.  Both proofs are combinatorial in nature.  In this note, we provide geometric proofs of both results, using affine Grassmannians.  We also  provide some additional results concerning minuscule coweights and the equidimensionality of the fibers of certain Bott-Samelson resolutions of affine Schubert varieties for ${\rm GL}_n$. 
\end{abstract}

\section{Introduction}
\markboth{Thomas J. Haines}{Structure constants for Hecke and representation rings}

Let $G$ denote a connected reductive group over an algebraic closure $k$ of a finite field ${\mathbb F}_q$.  We assume $G$ is defined and split over ${\mathbb F}_q$.  Let $T$ denote an ${\mathbb F}_q$-split maximal torus of $G$.   Let $B = TU$ denote an ${\mathbb F}_q$-rational Borel subgroup containing $T$.  Let $X_*(T)$ denote the set of cocharacters of $T$, and let $X^\vee_+$ denote the set of $B$-dominant elements of $X_*(T)$.  Let $\rho$ denote half the sum of the $B$-positive roots for $G$.  Let $\langle \cdot, \cdot \rangle : X^*(T) \times X_*(T) \rightarrow \ZZ$ denote the canonical pairing.

\medskip

Fix a prime $\ell$ different from the characteristic of $k$, and let $\hat{G} = \hat{G}(\bar{\QQ}_\ell)$ denote the Langlands dual of $G$ over the field $\bar{\QQ}_\ell$ (an algebraic closure of $\QQ_\ell$).  Each dominant cocharacter $\mu \in X^\vee_+$ can be thought of as a dominant weight for $\hat{G}$, hence such a cocharacter 
gives rise to a unique irreducible $\bar{\QQ}_\ell$-linear representation $V_\mu$ of 
$\hat{G}$ having highest weight $\mu$.  Let ${\rm Rep}(\hat{G})$ denote the category of rational representations of $\hat{G}$ over $\bar{\QQ}_\ell$.

\medskip

We will work with the spherical Hecke algebra of $G$ in the function-field setting.  Let ${\mathcal O} = k[[t]]$ and $F = k((t))$.  Also, let ${\mathcal O}_q = {\mathbb F}_q[[t]]$ and $F_q = {\mathbb F}_q((t))$.  
Let $K := G({\mathcal O})$ denote the ``maximal compact'' subgroup of $G(F)$.  Let $K_q = G({\mathcal O}_q)$.  We let ${\mathcal H}_q$ denote the Hecke algebra for $G$ over $F_q$, i.e., the convolution algebra of compactly-supported $K_q$-bi-invariant $\bar{\QQ}_\ell$-valued functions on $G(F_q)$ (convolution is defined using the Haar measure on $G(F_q)$ which assigns measure one to the compact open subgroup $K_q$).  For $\mu \in X^\vee_+$, define an element $f_\mu \in {\mathcal H}_q$ by
$$
f_\mu = {\rm char}(K_q \, t_\mu \, K_q)
$$
where $t_\mu := \mu(t) \in T(F_q)$.  By the Cartan decomposition, the elements $f_\mu$ form a $\bar{\QQ}_\ell$-basis for ${\mathcal H}_q$.  For later use let us recall that the convolution operation $*$ is defined, for $f_1,f_2 \in {\mathcal H}_q$,  by
$$
(f_1 * f_2)(g) = \int_{G(F_q)} f_1(x)f_2(x^{-1}g) dx.
$$

\medskip

Now let $\mu_1, \dots, \mu_r$ and $\lambda$ be dominant coweights for $G$.  Denote 
$\mu_\bullet = (\mu_1, \dots, \mu_r)$ and $|\mu_\bullet| = \mu_1 + \cdots + \mu_r$. 
We consider the decomposition in ${\rm Rep}(\hat{G})$
$$
V_{\mu_1} \otimes \cdots \otimes V_{\mu_r} = \bigoplus_{\lambda \in X^\vee_+ ,\,\, \lambda \leq |\mu_\bullet| } V_\lambda \otimes V^\lambda_{\mu_\bullet},
$$
where $V^\lambda_{\mu_\bullet}$ denotes the multiplicity vector space.  Here $\lambda_1 \leq \lambda_2$ means that $\lambda_2 - \lambda_1$ is a sum of $B$-positive coroots.  (It follows from the discussion in section 2 that if $V_\lambda$ occurs as a summand, then $\lambda \leq |\mu_\bullet|$ necessarily holds.)  The numbers ${\rm dim}(V^\lambda_{\mu_\bullet})$ are the {\em structure constants for the representation ring of $\hat{G}$}.

\medskip

In a parallel manner, for $\mu_\bullet$ and $\lambda$ as above, we write
$$
f_{\mu_1} * \cdots * f_{\mu_r} = \sum_\lambda c^\lambda_{\mu_\bullet} f_\lambda,
$$
for certain constants $c^\lambda_{\mu_\bullet}$; in this paper these are referred to 
as the {\em structure constants for the Hecke algebra ${\mathcal H}_q$}.  It is known and easy to see that the structure constants all belong to $\ZZ$.  In fact, as $q$ varies, they are given by polynomial functions in $q$ with integer coefficients: $c^\lambda_{\mu_\bullet} 
= c^\lambda_{\mu_\bullet}(q)$.  (See Lemmas 9.15, 9.18 of \cite{KLM}).

\medskip

The purpose of this paper is to study the relation between the following two properties of a collection $(\mu_\bullet, \lambda)$ of dominant cocharacters:

\medskip

${\bf Rep}(\mu_\bullet,\lambda)$: The irreducible representation $V_\lambda$ occurs with non-zero 
multiplicity in $V_{\mu_1} \otimes \cdots \otimes V_{\mu_r}$.  That is, ${\rm dim}(V^\lambda_{\mu_\bullet}) > 0$.

${\bf Hecke}(\mu_\bullet,\lambda)$: The function $f_\lambda$ appears with non-zero coefficient in $f_{\mu_1} * \cdots * f_{\mu_r}$.  That is, $c^\lambda_{\mu_\bullet} \neq 0$.  Equivalently, $K_q \, t_\lambda \, K_q \subset K_q \, t_{\mu_1} \, K_q \cdots K_q \, t_{\mu_r} \, K_q$.

\medskip

The following result was first proved for general groups $G$ by M. Kapovich, B. Leeb, and J. Millson in \cite{KLM}, Theorem 9.19.  In this paper we give another approach, based on the geometry of affine Grassmannians.  

\begin{theorem}[KLM] \label{Rep_to_Hecke}
If ${\bf Rep}(\mu_\bullet,\lambda)$ holds, then ${\bf Hecke}(\mu_\bullet,\lambda)$ also holds.  More precisely, we have the following description of the Hecke algebra structure constants
$$
c^{\lambda}_{\mu_\bullet}(q) = {\rm dim}(V^\lambda_{\mu_\bullet})\,  q^{\langle \rho, |\mu_\bullet| - \lambda \rangle} \, + \, \{ \mbox{terms with lower $q$-degree} \}.
$$

\end{theorem}

\medskip

For the general linear group, the converse implication also holds, and is originally due to P. Hall, using the combinatorics of Hall polynomials (cf. \cite{KLM}, section 9.6, and \cite{Mac}, part II, Theorem (4.3)).  We present a geometric proof here, again using affine Grassmannians.

\begin{theorem}[P. Hall] \label{Hecke_to_Rep}
Let $G = {\rm GL}_n$.  Then ${\bf Hecke}(\mu_\bullet, \lambda) \Rightarrow {\bf Rep}(\mu_\bullet,\lambda)$.
\end{theorem}

\medskip

The approach of this paper is to reformulate the two properties ${\bf Hecke}(\mu_\bullet,\lambda)$ and ${\bf Rep}(\mu_\bullet,\lambda)$ in terms of properties of the ``multiplication'' morphism
$$
m_{\mu_\bullet} : \widetilde{\mathcal Q}_{\mu_\bullet} \rightarrow \bar{\mathcal Q}_{|\mu_\bullet|}
$$
used to define the convolution of $K$-equivariant perverse sheaves on the affine Grassmannian.  Here the domain is a twisted product of the closures $\bar{\mathcal Q}_{\mu_i}$ of $K$-orbits ${\mathcal Q}_{\mu_i}$ in the affine Grassmannian for $G$.  The morphism $m_{\mu_\bullet}$ is given by forgetting all but the last factor in the twisted product.  It is a (stratified) semi-small, proper, and birational 
morphism (cf. section 2).  
One can view it as a partial desingularization of $\bar{\mathcal Q}_{|\mu_\bullet|}$, directly analogous to the Bott-Samelson partial desingularizations of Schubert varieties. (When the coweights $\mu_i$ are all minuscule, the domain is smooth, so in that case 
$m_{\mu_\bullet}$ is a genuine resolution of singularities.)   The semi-smallness 
of $m_{\mu_\bullet}$ means that the fiber over any point $y$ in a $K$-orbit stratum 
${\mathcal Q}_\lambda \subset \bar{\mathcal Q}_{|\mu_\bullet|}$ has dimension bounded above by $\langle \rho, |\mu_\bullet| - \lambda \rangle$.

\medskip

The two properties can be translated into properties of the fibers of $m_{\mu_\bullet}$.  The geometric Satake isomorphism (cf. \cite{Gi}, \cite{MV}, \cite{NP}) plays a key role in this step.

\medskip

\begin{prop}\label{summary} Suppose $(\mu_\bullet,\lambda)$ is such that ${\mathcal Q}_\lambda \subset \bar{\mathcal Q}_{|\mu_\bullet|}$ (equivalently, $\lambda \leq |\mu_\bullet|$).  Fix $y \in {\mathcal Q}_\lambda$.  Then:

1. $V^\lambda_{\mu_\bullet} \neq 0$ if and only if $m^{-1}_{\mu_\bullet}(y)$ has an irreducible 
component of dimension $\langle \rho, |\mu_\bullet| - \lambda \rangle$.  Moreover, such a component necessarily meets the open stratum ${\mathcal Q}_{\mu_\bullet}$.

2. Suppose $y$ is ${\mathbb F}_q$-rational.  Then $f_\lambda$ occurs in $f_{\mu_1} * \cdots * f_{\mu_n}$ if and only if $m^{-1}_{\mu_\bullet}(y)$ 
meets the open stratum ${\mathcal Q}_{\mu_\bullet}$.  In that case, we have 
$$
c^\lambda_{\mu_\bullet}(q) = \#(m_{\mu_\bullet}^{-1}(y) \cap {\mathcal Q}_{\mu_\bullet})({\mathbb F}_q) > 0.
$$
\end{prop}

\medskip

The first part of Theorem \ref{Rep_to_Hecke} follows immediately from Proposition \ref{summary}.  We deduce the second part concerning $c^\lambda_{\mu_\bullet}$ with little additional effort, by using the Weil conjectures to approximate the number of points on the algebraic varieties $m^{-1}_{\mu_\bullet}(y) \cap {\mathcal Q}_{\mu_\bullet}$.  

\medskip

We also have the following partial converse to Theorem \ref{Rep_to_Hecke}, valid for every group $G$.  
It is an easy and purely combinatorial consequence of the P-R-V conjecture (now a theorem due independently to S. Kumar and O. Mathieu, cf. \cite{Lit}).

\medskip

\begin{prop} \label{partial_converse}
If ${\bf Hecke}(\mu_\bullet,\lambda)$ holds, then there exist dominant coweights $\mu'_i \leq \mu_i$ ($1 \leq i \leq r$), such that ${\bf Rep}(\mu'_\bullet,\lambda)$ holds.
\end{prop}

\medskip

To prove Theorem \ref{Hecke_to_Rep}, we need finer information about the fibers of $m_{\mu_\bullet}$ in the ${\rm GL}_n$ case.  The first ingredient is the following proposition valid for all groups having minuscule coweights (recall that a coweight $\mu$ is {\em minuscule} provided that $\langle \alpha, \mu \rangle \in \{ -1, 0, 1 \}$, for every root $\alpha$ of $G$).  It is also a consequence of the P-R-V conjecture.

\medskip

 \begin{prop} \label{minuscule}
If $\mu_1, \dots, \mu_r$ are all minuscule coweights for $G$, and $\lambda \leq |\mu_\bullet|$, then both ${\bf Rep}(\mu_\bullet, \lambda)$ and ${\bf Hecke}(\mu_\bullet,\lambda)$ hold.
\end{prop}

\medskip

In view of Proposition \ref{summary}, this implies:

\begin{cor} \label{minuscule_relevance}
Let $\mu_1, \dots, \mu_r$ let be minuscule coweights for $G$, and let $\lambda \leq |\mu_\bullet|$.  Then for any $y \in {\mathcal Q}_\lambda$, we have 
$$
{\rm dim}(m_{\mu_\bullet}^{-1}(y)) = \langle \rho, |\mu_\bullet| - \lambda \rangle.
$$
In particular, every stratum of $\bar{\mathcal Q}_{|\mu_\bullet|}$ is relevant for the semi-small morphism $m_{\mu_\bullet}$.
\end{cor}

\medskip

We make essential use of the further information that in the ${\rm GL}_n$ case, the fibers above are {\em equidimensional}.

\medskip

\begin{prop} \label{equal_dim}
Let $G = {\rm GL}_n$.  Let $\mu_1, \dots, \mu_r$ be minuscule coweights.  Suppose ${\mathcal Q}_\lambda \subset \bar{\mathcal Q}_{|\mu_\bullet|}$, and let $y \in {\mathcal Q}_\lambda$.  Then every irreducible component of the fiber
$$
m^{-1}_{\mu_\bullet}(y)
$$
has dimension $\langle \rho, |\mu_\bullet| - \lambda \rangle$. 
\end{prop}

\medskip

The proof of Proposition \ref{equal_dim} proceeds by reduction to a theorem of N. Spaltenstein \cite{Sp} concerning partial Springer resolutions of the nilpotent cone for ${\rm GL}_n$.  The lack of an analogous result for other groups is one reason Proposition \ref{equal_dim} can be proved (at the moment) only for ${\rm GL}_n$.

\medskip

In fact for the general linear group, Proposition \ref{equal_dim} can be used to prove a seemingly stronger result.

\medskip

\begin{prop} \label{strong_equal_dim}
Let $\mu_1, \dots, \mu_r$ be dominant coweights for $G = {\rm GL}_n$.  Let $y \in {\mathcal Q}_\lambda \subset \bar{\mathcal Q}_{|\mu_\bullet|}$.  Let ${\mathcal Q}_{\mu'_\bullet} \subset \widetilde{\mathcal Q}_{\mu_\bullet}$ be the stratum indexed by $\mu'_\bullet = (\mu'_1, \dots, \mu'_r)$ for dominant coweights $\mu'_i \leq \mu_i$ ($1 \leq i \leq r$).  Then any irreducible component of the fiber $m^{-1}_{\mu_\bullet}(y)$ whose generic point belongs to ${\mathcal Q}_{\mu'_\bullet}$ has dimension $\langle \rho, |\mu'_\bullet| - \lambda \rangle$.
\end{prop}

\medskip

As we explain in section 8, Proposition \ref{strong_equal_dim} quickly implies Theorem \ref{Hecke_to_Rep}.

\medskip

Proposition \ref{strong_equal_dim} appears to be special to ${\rm GL}_n$.  Indeed, it is actually somewhat stronger than the implication ${\bf Hecke}(\mu_\bullet,\lambda) \Rightarrow {\bf Rep}(\mu_\bullet, \lambda)$, which is known to fail in general (e.g. for ${\rm SO}(5)$ or $G_2$, cf. \cite{KLM}, section 9.5).  But Proposition \ref{equal_dim} could remain valid if ${\rm GL}_n$ is replaced with an arbitrary group $G$, and it would be interesting to clarify the situation.   Such a generalization of Proposition 
\ref{equal_dim} would have applications to proving a type of ``Saturation theorem'' for a general reductive group (along the same lines as the new proof in \cite{KLM} of the Saturation theorem for ${\rm GL}_n$, originally proved by Knutson-Tao  \cite{KT}).

\medskip

\noindent {\em Remark}.  Suitably reformulated, the main results of this paper hold if the coefficient field $k = \bar{\mathbb F}_q$ is replaced by any algebraically closed field $\kappa$ (e.g. ${\mathbb C}$).  We let now ${\mathcal O} = \kappa [[t]]$, $K = G({\mathcal O})$, etc., and replace ${\bf Hecke}(\mu_\bullet,\lambda)$ with

\smallskip

${\bf Hecke'}(\mu_\bullet,\lambda): Kt_\lambda K \subset Kt_{\mu_1}K \cdots Kt_{\mu_r}K.$

\smallskip

\noindent Then the arguments of this paper prove that ${\bf Rep}(\mu_\bullet,\lambda) \Rightarrow 
{\bf Hecke'}(\mu_\bullet, \lambda)$ for every $G$, and that ${\bf Rep}(\mu_\bullet,\lambda) \Leftrightarrow {\bf Hecke'}(\mu_\bullet, \lambda)$ for ${\rm GL}_n$.  The statements concerning dimensions of fibers and their equidimensionality also remain valid.  Thus, here we avoid the hypothesis that $F_q$ have finite residue field, essential for Hall's proof of Theorem \ref{Hecke_to_Rep} using Hall polynomials, and for the proof of Theorem \ref{Rep_to_Hecke} in (\cite{KLM}, Theorem 9.19).

\bigskip

\noindent {\em Acknowledgments}.  It is a pleasure to thank John Millson for several stimulating conversations, and for delivering a series of lectures on \cite{KLM} at the University of Maryland which sparked my interest in this subject.   My intellectual debt to B.C. Ng\^{o} and P. Polo will be clear to the reader familiar with \cite{Ngo} or \cite{NP}, where several of the ideas in this paper appear, at least implicitly.  I also thank J.K. Yu for some helpful conversations, and M. Rapoport for his comments on the paper.  Finally, I am grateful to R. Kottwitz for suggesting some significant improvements in exposition, and for making detailed comments on the first version of this paper.

\section{Review of affine Grassmannians}

In this section we recall some well-known notions relating to affine Grassmannians.  The reader can find further details in \cite{Gi}, \cite{MV}, and \cite{NP}.

\medskip

We will work with the affine Grassmannian
$$
{\mathcal Q} = G(F)/K,
$$
which can be thought of as the $k$-points of an ind-scheme defined over 
${\mathbb F}_q$.  The group scheme $K = G({\mathcal O})$  acts naturally on 
${\mathcal Q}$ (on the left).  By the Cartan decomposition, the $K$-orbits are parametrized by dominant coweights $\lambda \in X^\vee_+$.   Indeed, let $e_0$ denote the base point of ${\mathcal Q}$ corresponding to the coset $K$, and let 
$e_\lambda = t_\lambda e_0$.  Then ${\mathcal Q}_\lambda = Ke_\lambda$ is the 
$K$-orbit corresponding to $\lambda$.  It is well-known that ${\mathcal Q}_\lambda$ is a smooth quasi-projective variety of dimension $\langle 2\rho, \lambda \rangle $, defined over ${\mathbb F}_q$.  Its closure $\bar{\mathcal Q}_\lambda$ in 
${\mathcal Q}$ is projective, but in general is not smooth.  

\medskip 
Let 
$j_\lambda : {\mathcal Q}_\lambda \hookrightarrow \bar{\mathcal Q}_\lambda$ denote the open immersion into the closure, and let us define the intersection complex
$$
{\mathcal A}_\lambda := {\rm IC}_\lambda(\bar{\QQ}_\ell) = 
j_{\lambda, !*}(\bar{\QQ}_\ell)[\langle 2 \rho, \lambda \rangle].
$$
Here we are applying $j_{\lambda,!*}$, the Goresky-MacPherson middle extension functor (for the middle perversity, cf. \cite{GM}, \cite{BBD}), to the shifted constant sheaf on the smooth variety ${\mathcal Q}_\lambda$.  The cohomological shift by the dimension of ${\mathcal Q}_\lambda$ is to ensure the result is a perverse sheaf.  It is known that ${\mathcal A}_\lambda$ is a self-dual $K$-equivariant simple perverse sheaf on ${\mathcal Q}$.

\medskip

Let $P_K({\mathcal Q})$ denote the category of $\bar{\QQ}_\ell$-linear $K$-equivariant perverse sheaves on ${\mathcal Q}$.  This is a semi-simple abelian category (cf. \cite{Ga}), whose simple objects are precisely the intersection complexes 
${\mathcal A}_\lambda$, for $\lambda \in X^\vee_+$.   In fact there exists a tensor (or ``convolution'') operation 
$$
\star : P_K({\mathcal Q}) \times P_K({\mathcal Q}) \rightarrow P_K({\mathcal Q})
$$
(defined below) which gives $P_K({\mathcal Q})$ the structure of a neutral Tannakian category over $\bar{\QQ}_\ell$.  The following theorem which identifies the corresponding algebraic group as the Langlands dual plays a key role in this paper.  We refer the reader to \cite{Gi},\cite{MV}, \cite{NP}, and the appendix of \cite{Nad}, for details of the proof.

\begin{theorem}[Geometric Satake Isomorphism] \label{geometric_Satake}
There is an equivalence of tensor categories
$$
(P_K({\mathcal Q}), \star) \,\, \widetilde{\longrightarrow} \,\, ({\rm Rep}(\hat{G}), \otimes),
$$
under which ${\mathcal A}_\lambda$ corresponds to the irreducible representation $V_\lambda$ with highest weight $\lambda$.
\end{theorem}

\medskip

\subsection{Definition of the twisted product}

To define the operation $\star$, we need some more preliminaries.
First, recall that any ordered pair of elements $L,L' \in {\mathcal Q}$ gives rise to an element ${\rm inv}(L,L') \in X^\vee_+$ (the ``relative position'' of $L,L'$).  The map ${\rm inv}: {\mathcal Q} \times {\mathcal Q} \rightarrow X^\vee_+$
 is defined as the composition
$$
{\rm inv} : G(F)/K \times G(F)/K \rightarrow G(F) \backslash \big[G(F)/K \times G(F)/K \big] = K \backslash G(F) /K \tilde{\rightarrow} X^\vee_+,
$$
where $G(F)$ acts diagonally on $G(F)/K \times G(F)/K$.
For $G = {\rm GL}_n$, ${\rm inv}(L,L')$ is just the usual relative position of two ${\mathcal O}$-lattices in $F^n$ given by the theory of elementary divisors.

\medskip

The usual partial ordering $\leq$ on the set of dominant coweights corresponds to the closure relation in ${\mathcal Q}$: ${\mathcal Q}_\lambda \subset \bar{\mathcal Q}_\mu$ if and only if $\lambda \leq \mu$.  It follows that ${\rm inv}(L,L') \leq \mu$ if and only if there exists $g \in G(F)$ such that $gL = e_0$ and $gL' \in \bar{\mathcal Q}_\mu$.  Thus, for $L \in {\mathcal Q}$ and $\mu \in X^\vee_+$ fixed, the set of $L'$ with ${\rm inv}(L,L') \leq \mu$ can be thought of as the $k$-points of a projective algebraic variety, isomorphic to $\bar{\mathcal Q}_\mu$.

\medskip

Now let $\mu_\bullet = (\mu_1, \dots, \mu_r)$,  where $\mu_i \in X^\vee_+$ for $1 \leq i \leq r$.  We define the twisted product scheme 
$$
\widetilde{\mathcal Q}_{\mu_\bullet} = \bar{\mathcal Q}_{\mu_1} \tilde{\times} \cdots 
\tilde{\times} \bar{\mathcal Q}_{\mu_n}
$$
to be the subscheme of ${\mathcal Q}^n$ consisting of points $(L_1, \dots, L_r)$ such that ${\rm inv}(L_{i-1},L_i) \leq \mu_i$ for $1 \leq i \leq r$ (letting $L_0 = e_0$).  Projection onto the last coordinate gives a proper surjective birational map
$$
m_{\mu_{\bullet}}: \widetilde{\mathcal Q}_{\mu_\bullet} \rightarrow \bar{\mathcal Q}_{|\mu_\bullet|}.
$$
The birationality follows from the more precise statement that the restriction of $m_{\mu_\bullet}$ to the inverse image of the open stratum ${\mathcal Q}_{|\mu_\bullet|} \subset \bar{\mathcal Q}_{|\mu_\bullet|}$ is an isomorphism.  This in turn can be deduced from the corresponding statement for an analogous {\em Demazure resolution} of an affine Schubert variety in the affine flag variety.  We omit the details.

\medskip

\subsection{Semi-small morphisms}

We shall make use of the notion of semi-small morphisms.  Let $f: X = \cup_\alpha X_\alpha \rightarrow Y = \cup_\beta Y_\beta$ be a proper surjective birational morphism between stratified spaces.  Suppose each $f(X_\alpha)$ is a union of strata $Y_\beta$.  We say $f$ is {\em semi-small} if, whenever $y \in Y_\beta \subset f(X_\alpha)$, then
$$
{\rm dim}(f^{-1}(y) \cap X_\alpha) \leq \frac{1}{2}({\rm dim}(X_\alpha) - {\rm dim}(Y_\beta)).$$
(In \cite{MV}, this notion is termed ``stratified semi-small''.  In the usual terminology (\cite{GM}), the domain of a semi-small morphism is assumed to be smooth.)

\medskip
We say $f$ is {\em locally trivial} (in the stratified sense) if whenever $Y_\beta \subset f(X_\alpha)$, the restriction of $f : f^{-1}(Y_\beta) \rightarrow Y_\beta$ to $X_\alpha \cap f^{-1}(Y_\beta)$ is Zariski-locally a trivial fibration.  In this case we have, for every $y \in Y_\beta \subset f(X_\alpha)$:
$$
{\rm dim}(f^{-1}(y) \cap X_\alpha) + {\rm dim}(Y_\beta) = {\rm dim}(f^{-1}(Y_\beta) \cap X_\alpha).
$$
 
\medskip

The target $\bar{\mathcal Q}_{|\mu_\bullet|}$ of $m_{\mu_\bullet}$ is stratified by the locally-closed subschemes ${\mathcal Q}_\lambda$  ($\lambda \leq |\mu_\bullet|$), and the 
domain $\bar{\mathcal Q}_{\mu_1} \tilde{\times} \cdots 
\tilde{\times} \bar{\mathcal Q}_{\mu_n}$ is stratified by the locally-closed subspaces 
${\mathcal Q}_{\mu'_\bullet} = {\mathcal Q}_{\mu'_1} \tilde{\times} \cdots 
\tilde{\times} {\mathcal Q}_{\mu'_n}$, where $\mu'_\bullet = (\mu'_1,\dots,\mu'_n)$ satisfies $\mu'_i \leq \mu_i$, for every $i$.  (The definition of the subspace is the same as that of the ambient space, except that the inequalities ${\rm inv}(L_{i-1},L_i) \leq \mu_i$ are replaced by the equalities ${\rm inv}(L_{i-1},L_i) = \mu'_i$.)  With respect to these stratifications, we have the following result.

\medskip

\begin{prop}[\cite{MV}, \cite{NP}]\label{semi-small}
The morphism $m_{\mu_\bullet}$ is a semi-small and locally trivial morphism.
\end{prop}

\medskip

\noindent {\em Remark.}  The proof of semi-smallness over ${\mathbb F}_q$ does not appear explicitly in the literature in precisely this generality.  The proof of Lemme 9.3 in \cite{NP} for the special case where all $\mu_i$ are minuscule or quasi-minuscule works as well for the general case, modulo the inequality ${\rm dim}(S_\nu \cap {\mathcal Q}_\lambda) \leq \langle \rho, \lambda + \nu \rangle$, where $\nu \in X_*(T)$ and $S_\nu := U(F)e_\lambda$.  In fact the equality ${\rm dim}(S_\nu \cap {\mathcal Q}_\lambda) = \langle \rho, \lambda + \nu \rangle$ can be deduced (as Ng\^{o} and Polo remark), from their Theoreme 3.1.  This indirect method to prove the dimension equality is not circular, since Ng\^{o} and Polo use the equality ${\rm dim}(S_\nu \cap {\mathcal Q}_\lambda) = \langle \rho, \lambda + \nu \rangle $ only in the special case mentioned above in their proof of Theoreme 3.1, and they prove this special case by direct means.

\medskip

\subsection{Convolution of perverse sheaves}

If ${\mathcal G}_i$ ($1 \leq i \leq r$) are elements of $P_K({\mathcal Q})$, choose coweights $\mu_i$ such that ${\rm supp}({\mathcal G}_i) \subset \bar{\mathcal Q}_{\mu_i}$, for every $i$ (technical aside: this is possible, as we can assume with no loss of generality that each ${\mathcal G}_i$ is supported on only one connected component of ${\mathcal Q}$).  There is a unique ``twisted product'' perverse sheaf ${\mathcal G}_1 \tilde{\boxtimes} \cdots \tilde{\boxtimes} {\mathcal G}_r$ on the twisted product $\tilde{Q}_{\mu_\bullet}$, which is locally isomorphic to the usual exterior product ${\mathcal G}_1 \boxtimes \cdots \boxtimes {\mathcal G}_r$ on the product space $\bar{\mathcal Q}_{\mu_1} \times \cdots \times \bar{\mathcal Q}_{\mu_r}$.  (See 7.4 of \cite{HN}, or \cite{NP}, section 2, for another construction).  We then define
$$
{\mathcal G}_1 \star \cdots \star {\mathcal G}_r = Rm_{\mu_\bullet,*}({\mathcal G}_1 \tilde{\boxtimes} \cdots \tilde{\boxtimes} {\mathcal G}_r).
$$
This belongs to $P_K({\mathcal Q})$ by the semi-smallness of $m_{\mu_\bullet}$, and is independent of the choice of $\mu_\bullet$.

\medskip 

The important lemma below follows directly from the definitions and the following characterization of the intersection complex (\cite{BBD}, 2.1.17):   Let $X^0$ denote the open (smooth) stratum in a stratified variety $X = \cup_\alpha X_\alpha$.  Then ${\rm IC}(X)$ is the unique self-dual perverse extension of $\bar{\mathbb Q}_\ell[{\rm dim}(X)]$ on $X^0$ which satisfies, for each stratum $X_\alpha \neq X^0$, the property ${\mathcal H}^i{\rm IC}(X)|_{X_\alpha} = 0$ for $i \geq -{\rm dim}(X_\alpha)$.   

\medskip

\begin{lemma}\label{IC_upstairs}
${\mathcal A}_{\mu_1} \star \cdots \star {\mathcal A}_{\mu_r} = Rm_{\mu_\bullet, *} 
{\rm IC}(\widetilde{\mathcal Q}_{\mu_\bullet})$.
\end{lemma}

\medskip

We are going to use this to give a geometric description of the multiplicity space 
$V^\lambda_{\mu_\bullet}$, which will be the key step to proving Proposition \ref{summary}.

\section{Proof of Proposition \ref{summary}, Part (1)}

The key result is the following well-known proposition.  It has appeared without proof in several published articles (e.g. \cite{NP}).  We provide a proof here for the convenience of the reader.

\medskip

\begin{prop}\label{key}
Fix $\lambda \leq |\mu_\bullet|$, and $y \in {\mathcal Q}_\lambda$.  Then the space 
$V^\lambda_{\mu_\bullet}$ has a basis in canonical bijection with the set of irreducible components of $m^{-1}_{\mu_\bullet}(y)$ having maximal possible dimension, that is, $\langle \rho, |\mu_\bullet| - \lambda \rangle$.
\end{prop}

\medskip

Let us assume Proposition \ref{key} and deduce Proposition \ref{summary}, Part (1).
The first statement is obvious, so we verify the second.  But it is clear that an irreducible component $C$ of $m_{\mu_\bullet}^{-1}(y)$ with dimension $\langle \rho, |\mu_\bullet| - \lambda \rangle$ meets the open stratum ${\mathcal Q}_{\mu_\bullet}$ (if not, then $C$ meets only strictly smaller strata ${\mathcal Q}_{\mu'_\bullet}$, and then by Proposition \ref{semi-small}, the dimension of $C$ would be bounded above by
$$
{\rm dim}(m^{-1}_{\mu_\bullet}(y) \cap {\mathcal Q}_{\mu'_\bullet}) \leq \langle \rho, |\mu'_\bullet| - \lambda \rangle < \langle \rho, |\mu_\bullet| - \lambda \rangle.)
$$

\medskip

\noindent {\em Proof of Proposition \ref{key}:}
By Theorem \ref{geometric_Satake} we have 
$$
{\mathcal A}_{\mu_1} \star \cdots \star {\mathcal A}_{\mu_r} = \bigoplus_{\lambda \leq |\mu_\bullet|} {\mathcal A}_\lambda \otimes V^\lambda_{\mu_\bullet}.
$$

\medskip

Let $d = - {\rm dim}({\mathcal Q}_\lambda)$.  Now apply ${\mathcal H}^d_y$ (take the $d$th cohomology stalk at $y$) to both sides of this equation.

\medskip

The right hand side gives the vector space $V^\lambda_{\mu_\bullet}$, since ${\mathcal A}_\lambda$ restricted to ${\mathcal Q}_\lambda$ is $\bar{\QQ}_{\ell}[{\rm dim}({\mathcal Q}_\lambda)]$ (the constant sheaf placed in degree $-{\rm dim}({\mathcal Q}_\lambda)$) and since by definition of intersection complexes, ${\mathcal H}^d_y{\mathcal A}_{\lambda'}$ is zero if $\lambda < \lambda'$.  (In general, 
${\mathcal H}^i{\rm IC(X)}|_{X_{\alpha}} = 0$ if $i > - {\rm dim(X_\alpha)}$ and also for $i \geq - {\rm dim}(X_\alpha)$ if $X_\alpha$ is not the open stratum in $X$.)

\medskip

By Lemma \ref{IC_upstairs}, on the left hand side we get 
$$
H^{d}(m^{-1}_{\mu_\bullet}(y), {\rm IC}(\widetilde{\mathcal Q}_{\mu_\bullet})).
$$

\medskip

The result now follows by virtue of the following lemma.

\medskip

\begin{lemma}
Let $f: X = \cup_\alpha X_\alpha \rightarrow Y = \cup_\beta Y_\beta$ be a semi-small morphism between proper stratified schemes, and suppose $y \in Y_\beta$.  Let $d = - {\rm dim}(Y_\beta)$.  Then there is a canonical isomorphism
$$
H^d(f^{-1}(y),{\rm IC}(X)) = \bar{\mathbb Q}_\ell^{C_{max}(y)},
$$
where $C_{max}(y)$ is the set of irreducible components of $f^{-1}(y)$ having the maximal possible dimension $\frac{1}{2}({\rm dim}(X) - {\rm dim}(Y_\beta))$.
\end{lemma}

\begin{proof}

Let $X^0$ denote the open (smooth) stratum of $X$; recall that ${\rm IC}(X)|_{X^0} = \bar{\mathbb Q}_\ell[{\rm dim}(X)]$.  

Suppose we fix $X_\alpha \neq X^0$.  We claim that
$$
H^d(f^{-1}(y) \cap X_\alpha, {\rm IC}(X)) = 0.
$$
Indeed, if not, then the ``local-global'' spectral sequence
$$
H^p(f^{-1}(y) \cap X_\alpha, {\mathcal H}^q{\rm IC}(X)) \Rightarrow H^{p+q}(f^{-1}(y) \cap X_\alpha, {\rm IC}(X))$$
shows that there exists $q$ and $p=d-q$ such that the initial term is non-zero.  We have therefore $q < - {\rm dim}(X_\alpha)$, which together with the semi-smallness of $f$ 
gives
$$
p \leq 2 ({\rm dim}(f^{-1}(y) \cap X_\alpha)) \leq {\rm dim}(X_\alpha) - {\rm dim}(Y_\beta) < -q - {\rm dim}(Y_\beta),
$$
which is clearly impossible since $p + q = d = -{\rm dim}(Y_\beta)$.

\medskip

The same argument shows that $H^d(f^{-1}(y)\cap (X \backslash X^0), {\rm IC}(X)) = 0$.

\medskip

Since $f^{-1}(y)$ and $f^{-1}(y) \cap (X \backslash X^0)$ are proper, their cohomology and compactly-supported cohomology agree.  It follows from the above remarks that
$$
H^d(f^{-1}(y),{\rm IC}(X)) = H^d_c(f^{-1}(y) \cap X^0, \bar{\mathbb Q}_\ell[{\rm dim}(X)]) = 
H_c^{d + {\rm dim}(X)}( f^{-1}(y) \cap X^0, \bar{\mathbb Q}_\ell),
$$
from which the lemma follows easily.
\end{proof}

\section{Proof of Proposition \ref{summary}, Part (2)}

The next lemma follows from the definition of convolution in ${\mathcal H}_q$, and is left to the reader.

\medskip

\begin{lemma} \label{fiber_cardinality}
$c^\lambda_{\mu_\bullet}(q) = \#(m_{\mu_\bullet}^{-1}(y) \cap {\mathcal Q}_{\mu_\bullet})({\mathbb F}_q)$.  
\end{lemma}

\medskip

Therefore, if $f_\lambda$ occurs in $f_{\mu_1} * \cdots * f_{\mu_r}$, we obviously have 
$m_{\mu_\bullet}^{-1}(y) \cap {\mathcal Q}_{\mu_\bullet} \neq \emptyset$.

\medskip 

Conversely, suppose that $m_{\mu_\bullet}^{-1}(y) \cap {\mathcal Q}_{\mu_\bullet} \neq \emptyset$.  It follows from this that $c^\lambda_{\mu_\bullet}(q^n) > 0$ for $n >\!>0$. 
Using the definition of convolution, this means that
$$
t_\lambda \in K_{q^n} t_{\mu_1} K_{q^n} \cdots K_{q^n} t_{\mu_r} K_{q^n}.
$$
But by the discussion below, this condition can be expressed purely in terms of the extended affine Weyl group, and in particular, it is independent of $n$.  Therefore, it holds for $n=1$ and thus 
$c^\lambda_{\mu_\bullet}(q) > 0$. 

\medskip

To complete the argument proving independence of $n$, let $\widetilde{W} = X_*(T) \rtimes W$ denote the extended affine Weyl group.  Denote the translation element in $\widetilde{W}$ corresponding to $\mu \in X_*(T)$ also by the symbol $t_\mu$.  Let $I_{q^n} \subset K_{q^n}$ be the Iwahori subgroup defined to be the inverse image of $B$ under the homomorphism $G({\mathbb F}_{q^n}[[t]]) \rightarrow G({\mathbb F}_{q^n})$ given by $t \mapsto 0$.  For simplicity, write $K = K_{q^n}$ and $I = I_{q^n}$.

\medskip

Then $K = IWI := 
\coprod_{w \in W} IwI$ and $G({\mathbb F}_{q^n}((t))) = I\widetilde{W}I := \coprod_{w \in \widetilde{W}} IwI$.  Furthermore, standard results for BN pairs yield the identity
$$
IxIyI \subset \coprod_{\tilde{y} \preceq y} Ix\tilde{y}I,
$$
for $x,y \in \widetilde{W}$, where $\tilde{y}$ ranges over elements preceding $y$ in the Bruhat order on $\widetilde{W}$.  Using this  
we see 
\begin{align*}
t_\lambda \in Kt_{\mu_1}K \cdots Kt_{\mu_r}K &\Leftrightarrow t_\lambda \in IWt_{\mu_1}WI \cdots IWt_{\mu_r}WI \\
&\Leftrightarrow \exists \, x_i \in Wt_{\mu_i}W \,\, \mbox{such that} \,\, t_\lambda \in Ix_1I \cdots Ix_rI.
\end{align*}
But by standard facts for BN-pairs, the set ${\mathcal S} \subset \widetilde{W}$ appearing in the union
$$
Ix_1I \cdots Ix_rI = \coprod_{w \in {\mathcal S}} IwI
$$
depends only on the elements $x_1, \dots, x_r \in \widetilde{W}$ (and not on the power $q^n$ in $I = I_{q^n}$).

\section{End of proof of Theorem \ref{Rep_to_Hecke}}

It remains to prove the formula
$$
c^\lambda_{\mu_\bullet}(q) = {\rm dim}(V^\lambda_{\mu_\bullet}) \, q^{\langle \rho, |\mu_\bullet| - \lambda \rangle} + \{ \mbox{terms of lower $q$-degree} \}.
$$

Clearly we may prove this after base extension (enlarging $q$), so that we may assume the irreducible components of $m_{\mu_\bullet}^{-1}(y) \cap {\mathcal Q}_{\mu_\bullet}$ are defined over ${\mathbb F}_q$.  But then taking Lemma \ref{fiber_cardinality} and Proposition \ref{key} into account, the formula follows immediately from the following lemma, itself a consequence of Deligne's paper \cite{Weil2}.

\medskip

\begin{lemma}\label{growth}
If $X$ is a geometrically irreducible ${\mathbb F}_q$-variety, then the function 
$\#X({\mathbb F}_q)$ is of the form $q^{{\rm dim}(X)} + r(q)$, where 
$|r(q)| \leq O(q^{{\rm dim}(X) - 1/2})$. 
\end{lemma}

\section{The P-R-V conjecture and Propositions \ref{partial_converse} and \ref{minuscule}}

The Parthasarathy-Ranga-Rao-Varadarajan (P-R-V) conjecture has been proved independently by S. Kumar and O. Mathieu.  See \cite{Lit}, section 10, for a short proof using the Littelmann path model.

\medskip

\begin{theorem}[P-R-V Conjecture]\label{PRV}
For $1 \leq i \leq r$, let $\mu_i$ be a dominant coweight, and let $w_i \in W$ be an element of the finite Weyl group.  Suppose $\nu = w_1\mu_1 + \cdots + w_r \mu_r$ is dominant.  Then $V_\nu$ occurs with multiplicity at least one in the tensor product
$$
V_{\mu_1} \otimes \cdots \otimes V_{\mu_r}.
$$
\end{theorem}

\medskip
(This is usually stated only in the case $r=2$, but the above version follows easily by induction on $r$.)

\medskip

Theorem \ref{PRV} is the main ingredient to proving Proposition \ref{partial_converse} and the assertion ${\bf Rep}(\mu_\bullet,\lambda)$ in Proposition \ref{minuscule}.  Indeed, for the latter we only need to verify the following lemma to be able to take $\nu = \lambda$ in Theorem \ref{PRV}.

\medskip

\begin{lemma}\label{tensor_for_minuscule}
Suppose that each dominant coweight $\mu_i$ is minuscule, and that $\lambda$ is dominant and satisfies $\lambda \leq \mu_1 + \cdots + \mu_r$.  Then 
there exist elements $w_i \in W$ ($1 \leq i \leq r$), such that
$$
\lambda = w_1\mu_1 + \cdots + w_r\mu_r.
$$
\end{lemma}

\medskip

For $V \in {\rm Rep}(\hat{G})$, let $\Omega(V) \subset X^*(\hat{T}) = X_*(T)$ denote the set of its weights with respect to the dual torus $\hat{T}$.  Recall that 
$$
\Omega(V_\mu) = \{ \nu \in X^*(\hat{T}) ~ | ~ w\nu \leq \mu, \,\, \forall w \in W \}.
$$
The following well-known lemma is left to the reader.

\medskip

\begin{lemma} \label{weights}
Let $\mu_1, \dots, \mu_r$ be dominant coweights (also viewed as weights for $\hat{G}$).  Then 
$$
\Omega(V_{|\mu_\bullet|}) = \Omega(V_{\mu_1} \otimes \cdots \otimes V_{\mu_r}) = \{ \nu_1 + \cdots + \nu_r ~ | ~ \nu_i \in \Omega(V_{\mu_i}), \, 1 \leq i \leq r \}.
$$

\end{lemma}

\medskip
Note that any dominant $\lambda$ satisfying $\lambda \leq |\mu_\bullet|$ is necessarily a weight for $V_{|\mu_\bullet|}$.  Also, $\mu$ minuscule implies that $\Omega(V_\mu) = W\mu$.  Thus, Lemma \ref{weights} implies Lemma \ref{tensor_for_minuscule}.

\medskip

In the same way, we get the following result, proving  Proposition \ref{partial_converse}.

\medskip

\begin{prop} \label{application_of_PRV}
Suppose $\lambda \leq |\mu_\bullet|$ is dominant.  Then there exist dominant coweights $\mu'_i \leq \mu_i$ and elements $w_i \in W$ ($1 \leq i \leq r$) such that $\lambda = w_1\mu'_1 + \cdots + w_r\mu'_r$.  Consequently, the P-R-V conjecture implies that 
${\bf Rep}(\mu'_\bullet,\lambda)$ holds.  In particular, we have
$$
{\bf Hecke}(\mu_\bullet,\lambda) \Longrightarrow {\bf Rep}(\mu'_\bullet,\lambda).
$$
\end{prop}

\medskip

As for the assertion ${\bf Hecke}(\mu_\bullet,\lambda)$ of Proposition \ref{minuscule}, at this point one could appeal to Theorem \ref{Rep_to_Hecke}, but we prefer to give a direct proof.  Indeed, the assertion is obvious from the remark that $m_{\mu_\bullet}({\mathcal Q}_{\mu_\bullet}) = \bar{\mathcal Q}_{|\mu_\bullet|}$.  This in turn follows from the surjectivity of $m_{\mu_\bullet}$, since $\mu_i$ minuscule for every $i$ means $\bar{\mathcal Q}_{\mu_i} = {\mathcal Q}_{\mu_i}$ and thus $\widetilde{\mathcal Q}_{\mu_\bullet} = 
{\mathcal Q}_{\mu_\bullet}$.

\medskip

\section{Spaltenstein-Springer varieties and Proposition \ref{equal_dim}}

We begin with the statement of a crucial theorem of Spaltenstein \cite{Sp}.  
Let $V$ denote a $k$-vector space of dimension $d$, and let $(d_1, \dots, d_r)$ denote an ordered $r$-tuple of non-negative integers such that $d_1 + \cdots + d_r = d$.  The $r$-tuple $d_\bullet$ determines a standard parabolic subgroup $P \subset {\rm GL}(V)$.  Let us consider the variety of partial flags of type $P$:
$$
{\mathcal P} = \{ V_\bullet = (V = V_0 \supset V_1 \supset \cdots \supset V_r = 0) ~ | ~ 
{\rm dim}(V_{i-1}/V_i) = d_i, \,\, 1 \leq i \leq r \},
$$
which is isomorphic to ${\rm GL}(V)/P$. 

\medskip

For any nilpotent endomorphism $T \in {\rm End}(V)$, let ${\mathcal P}^T$ denote the closed subvariety of ${\mathcal P}$ consisting of partial flags $V_\bullet$ such that $T$ stabilizes each $V_i$.  This is simply the fiber over $T$ of the partial Springer resolution
$$
\pi: \widetilde{N}^P \rightarrow {\mathcal N},
$$
where ${\mathcal N} \subset {\rm End}(V)$ is the nilpotent cone, $\widetilde{\mathcal N}^P = \{ (T,V_\bullet) \in {\mathcal N} \times {\mathcal P} ~ | ~ V_\bullet \in {\mathcal P}^T \}$, and the morphism $\pi$ is the obvious forgetful one.

\medskip

Inside ${\mathcal P}^T$ we may consider the closed subvariety ${\mathcal P}^T_{\rm min}$ consisting of partial flags $V_\bullet \in {\mathcal P}^T$ such that 
$$
T \, \mbox{ acts by zero on} \,  V_{i-1}/V_i \,\, (1 \leq i \leq r).
$$
We call such varieties ${\mathcal P}^T_{\rm min}$ the {\em Spaltenstein-Springer varieties}.  
\medskip

We have the following fundamental result of N. Spaltenstein (\cite{Sp}, final Corollary).

\medskip

\begin{theorem}[Spaltenstein] \label{Spaltenstein}
The irreducible components of a Spaltenstein-Springer variety ${\mathcal P}^T_{\rm min}$ all have the same dimension.
\end{theorem}

\medskip

Now we turn to the proof of Proposition \ref{equal_dim}, which we will reduce to Theorem \ref{Spaltenstein}.

\medskip

Recall that $\mu_\bullet = (\mu_1, \dots, \mu_r)$, where each $\mu_i$ is a dominant minuscule coweight for ${\rm GL}_n$.  Working in the affine Grassmannian for ${\rm GL}_n$ and fixing $y \in {\mathcal Q}_\lambda \subset \bar{\mathcal Q}_{|\mu_\bullet|}$, we want to show that all the irreducible 
components of $m_{\mu_\bullet}^{-1}(y)$ have the same dimension.  Without loss of generality, we may assume $\mu_i = (1^{d_i},0^{n-d_i})$, where $1 \leq d_i \leq n-1$.  Write $d:= d_1 + \cdots + d_r$.  Suppose the point $y \in {\mathcal Q}_\lambda$ corresponds to the ${\mathcal O}$-lattice $L \subset F^n$.  We have $L \subset {\mathcal O}^n$.  Let $T$ denote the nilpotent endomorphism of the $d$-dimensional $k$-vector space $V:= {\mathcal O}^n/L$ induced by multiplication by $t$.  Then the fiber $m_{\mu_\bullet}^{-1}(y)$ consists of chains of $k$-vector spaces 
$$
{\mathcal O}^n = L_0 \supset L_1 \supset \cdots \supset L_r = L
$$
satisfying 
\begin{itemize}
\item ${\rm dim}(L_{i-1}/L_i) = d_i, \,\, \mbox{for each $i$}$,
\item $T \,\, \mbox{preserves each $L_i$}$,
\item $T \,\, \mbox{acts by zero on each $L_{i-1}/L_i$}$.
\end{itemize}
In other words, $m_{\mu_\bullet}^{-1}(y)$ is a Spaltenstein-Springer variety for the $d$-dimensional $k$-vector space $V$ and the nilpotent operator $T$.  Thus, Proposition \ref{equal_dim} follows from Theorem \ref{Spaltenstein}.
 
\medskip

\section{Proofs of Proposition \ref{strong_equal_dim} and Theorem \ref{Hecke_to_Rep}}

\noindent{\em Proof of Proposition \ref{strong_equal_dim}}:  We need the following set-up in order to bring Proposition \ref{equal_dim} into play.  Write each coweight $\mu_i$ as a sum of dominant minuscule coweights
$$
\mu_i = \nu_{i1} + \nu_{i2} + \cdots
$$
We consider the following diagram
$$
\xymatrix{
({\mathcal Q}_{\nu_{11}} \tilde{\times} \cdots ) \tilde{\times} \cdots \tilde{\times} 
({\mathcal Q}_{\nu_{r1}} \tilde{\times} \cdots ) \ar[d]^{\eta} \\
\bar{\mathcal Q}_{\mu_1} \tilde{\times} \cdots \tilde{\times} \bar{\mathcal Q}_{\mu_r} \ar[d]^{m_{\mu_\bullet}} \\
\bar{\mathcal Q}_{|\mu_\bullet|},
}
$$
where $\eta = m_{\nu_{1\bullet}} \tilde{\times} \cdots \tilde{\times} m_{\nu_{r\bullet}}$.  The composition $m_{\nu_{\bullet \bullet}} = m_{\mu_\bullet} \circ \eta$ is just the usual forgetful morphism, once the domain is identified with the twisted product of the varieties ${\mathcal Q}_{\nu_{ij}}$.  In particular, all three morphisms $m_{\mu_\bullet}$, $\eta$, and $m_{\nu_{\bullet \bullet}}$ are semi-small, proper, birational, surjective, and locally trivial in the stratified sense.  Moreover, the conclusion of Proposition \ref{equal_dim} applies to $m_{\nu_{\bullet \bullet}}$.  Note that $|\nu_{\bullet \bullet}| = |\mu_\bullet|$.

\medskip
Let $y \in {\mathcal Q}_\lambda \subset \bar{\mathcal Q}_{\mu_\bullet}$.  
Suppose $C'$ is an irreducible component of $m_{\mu_\bullet}^{-1}(y)$ such that $C' \cap {\mathcal Q}_{\mu'_\bullet}$ is an open dense subset in $C'$, so that ${\rm dim}(C' \cap {\mathcal Q}_{\mu'_\bullet}) = {\rm dim}(C')$.  Let $C$ be any irreducible component of $\eta^{-1}(C')$ which dominates $C'$.  Then $C \cap \eta^{-1}({\mathcal Q}_{\mu'_\bullet})$ is non-empty, hence is an open dense subset of $C$.  
Thus, using the semi-smallness and local triviality of $\eta$, we get
\begin{align*}
{\rm dim}(C) &= {\rm dim}(C \cap \eta^{-1}({\mathcal Q}_{\mu'_\bullet})) \\ 
             &\leq {\rm dim}(\eta^{-1}(C' \cap {\mathcal Q}_{\mu'_\bullet})) \\
             &\leq \langle \rho, |\nu_{\bullet \bullet}| - |\mu'_{\bullet}| \rangle + {\rm dim}(C' \cap {\mathcal Q}_{\mu'_\bullet}) \\
             &= \langle \rho, |\mu_\bullet| - |\mu'_\bullet| \rangle + {\rm dim}(C').
\end{align*}

Now $C$ is an irreducible component of $m_{\nu_{\bullet \bullet}}^{-1}(y)$, hence by Proposition \ref{equal_dim} has dimension $\langle \rho, |\mu_\bullet| - \lambda \rangle$.  Therefore 
$$
{\rm dim}(C') \geq \langle \rho, |\mu'_\bullet| - \lambda \rangle,
$$
and the desired equality follows from the semi-smallness of $m_{\mu_\bullet}$.  

\medskip

\noindent {\em Proof of Theorem \ref{Hecke_to_Rep}}:

This follows immediately from Proposition \ref{summary}, Proposition \ref{strong_equal_dim} and the remark that since 
${\mathcal Q}_{\mu_\bullet}$ is open in $\widetilde{\mathcal Q}_{\mu_\bullet}$, any irreducible component $C$ of the fiber $m_{\mu_\bullet}^{-1}(y)$ which meets ${\mathcal Q}_{\mu_\bullet}$ has $C \cap {\mathcal Q}_{\mu_\bullet}$ as a dense open subset.  This completes the proof of Theorem \ref{Hecke_to_Rep}.

\obeylines
Mathematics Department
University of Maryland
College Park, MD 20742-4015
USA

tjh@math.umd.edu

\end{document}